\documentclass[12pt,a4paper]{amsart}
\usepackage{amssymb,amsmath} %,amsmath,mathrsfs,amsfonts,amsthm
\usepackage[colorlinks, linkcolor=blue,  citecolor=blue,urlcolor=blue]{hyperref}
\usepackage{enumerate}
\newtheorem{theorem}{Theorem}[section] % 1st argument is your name for it
\newtheorem{lemma}[theorem]{Lemma}     % 2nd argument is what is printed

\newtheorem{remark}[theorem]{Remark}
\title[Explicit Barenblatt Profiles]{Explicit Barenblatt Profiles for Fractional Porous Medium Equations} % This is the full title of the paper
\author{Yanghong Huang}
\thanks{Department of Mathematics, Imperial College London, London SW7 2AZ, United Kingdom. 
 Email: yanghong.huang@imperial.ac.uk}
\begin{document}
\maketitle

\begin{abstract}
Several one-parameter families of explicit self-similar solutions are constructed for the porous medium equations with fractional operators. 
The corresponding self-similar profiles, also called \emph{Barenblatt profiles}, have the same forms as those of the classic porous medium 
equations. These new exact solutions complement current theoretical analysis of the underlying equations and are expected to 
provide insights for further quantitative investigations.
\end{abstract}

%\part{Use this type of header for very long papers only}
% use lowercase except for proper names

\section{Introduction}
The realistic modelling of phenomena in nature and science is usually described by nonlinear Partial Different Equations (PDEs). Comparing to their simplified linear counterparts, these nonlinear PDEs in general has no explicit representation of the solutions in terms of initial  and/or boundary conditions. Special explicit solutions, if available, are
often associated to certain symmetry groups of the underlying equation~\cite{bluman2002symmetry,olver2000applications}, 
including the most important ones, the scaling symmetry induced self-similar solutions.

Although self-similar solutions arise as exact solutions only with compatible initial and boundary conditions, they possess a unique position in the general theory of nonlinear partial differential equations. Take the Porous Medium Equation (PME) $u_t = \Delta u^m$ in $\mathbb{R}^N$ for example. As summarized in the monographs~\cite{vazquez2006smoothing,vazquez2006porous}, 
the self-similar solutions, also called Barenblatt-Kompaneets-Pattle-Zel'dovich solutions, characterize the long time asymptotic behaviours with nonnegative initial data; they indicate the parameter regimes  where the finite versus infinite speed of propagation of information is expected; they also provide a guidance to more refined questions like optimal regularity and optimal constants in various functional identities and inequalities. 

In this paper, we investigate the existence of certain explicit self-similar solutions of the porous medium equations with fractional operators, i.e.,
\begin{subequations}\label{eq:fpme}
    \begin{equation}\label{eq:fpme1}
        u_t + (-\Delta)^{s} u^m = 0,
\end{equation}
and
\begin{equation}\label{eq:fpme3}
    u_t =\nabla \cdot \big(u^{m-1}\nabla (-\Delta)^{-s} u\big).
\end{equation}
\end{subequations} 
The definition and related properties of the fractional Laplacian $(-\Delta)^s$ and its inverse
$(-\Delta)^{-s}$, together with the associated function spaces, can be found in the monographs~\cite{MR0350027,MR0290095} or the survey paper~\cite{MR2944369}.
When $s=2$ in~\eqref{eq:fpme1} or $s=0$ in~\eqref{eq:fpme3}, the classical PME is recovered (with different diffusion coefficients). The latter is also closely related to 
another variant with fractional  pressure
\begin{equation}\label{eq:fpme2}
    u_t = \nabla\cdot\big( u\nabla (-\Delta)^{-s} u^{m-1}\big).
\end{equation}
In fact,~\eqref{eq:fpme3} coincides with~\eqref{eq:fpme2}, when $m=2$ in both cases. 

Despite the equivalence of~\eqref{eq:fpme1},~\eqref{eq:fpme3} and~\eqref{eq:fpme2} 
to the classical PME in some ranges of $s$ and $m$, the three equations exhibit quite different qualitative properties. The basic theory of~\eqref{eq:fpme1} is studied in~\cite{MR2737788} for $s=1/2$ and in~\cite{de2012general} for general $s \in (0,1)$, followed by more refined quantitative estimates~\cite{bonforte2012quantitative,vazquez2013optimal,JLVBV}. In contrast, the notable feature of~\eqref{eq:fpme2} is the finite speed of propagation, studied for $m=2$ by Caffarelli and V{\'a}zquez~\cite{MR2773189,caffarelli2011nonlinear} and 
for general $m>1$ by Biler, Imbert and Karch~\cite{biler2011barenblatt,biler2013nonlocal}. The  variant~\eqref{eq:fpme3} has been studied only recently~\cite{prop}; depending on $m$, the equation can have both finite (for $1<m<2$) and infinite speed of propagation (for $m>2$).

One of the most important approaches to the study of qualitative and quantitative properties of PDEs
is to examine their self-similar solutions, whenever they exist. 
The self-similar solutions are related to the scaling symmetry groups of the PDEs, leading to transformed equations in scale-invariant similarity variables. After the reduction using similarity variables, the resulting equations for the self-similar profiles, called \emph{Barenblatt profiles} below,  still inherit some of the remaining scaling symmetries. As summarized in~\cite{MR1426127}, for self-similarity of the first kind, the scaling exponents can be determined a priori and explicit Barenblatt profiles can often be obtained. For self-similarity of the \emph{second kind}, also called \emph{anomalous scaling}, Barenblatt profiles are in general not available, because of the unknown anomalous exponents. Second kind self-similarity can be demonstrated by the PME $u_t=\Delta u^m$ in the fast diffusion regime $m<m_c=(N-2)_+/N$ where solutions
are known to vanish in finite time. Although no explicit Barenblatt profiles are expected in this case, the remaining scaling symmetry of the reduced equation allows one to give a detailed phase plane analysis to study the existence, uniqueness and monotonicity of the  profiles~\cite{MR1348964,vazquez2006smoothing}.

Unfortunately, there is limited usage of the remaining scaling symmetry of profile equations from the fractional porous medium equations~\eqref{eq:fpme}, for both first and second kind self-similarities, because the local characterization of Lie symmetry~\cite{bluman2002symmetry,olver2000applications}   
is destroyed by the nonlocal operator. As a consequence, explicit Barenblatt profiles are much more difficult to find. Surprisingly, all Barenblatt profiles of~\eqref{eq:fpme2} are obtained by Biler, Imbert and Karch~\cite{biler2011barenblatt,biler2013nonlocal} for any $s\in (0,1)$ and $m>1$, which are shown to be proportional to $(R^2-|y|^2)_+^{\frac{1-s}{m-1}}$ for some $R>0$. 
In this paper, we will focus on the less known explicit profiles for~\eqref{eq:fpme1} and~\eqref{eq:fpme3}, despite the existence, uniqueness and many qualitative properties presented in~\cite{vazquez2012barenblatt}
for~\eqref{eq:fpme1}.
In contrast to the explicit two-parameter family (for any $s$ and $m$) of profiles for~\eqref{eq:fpme2}, we can only find isolated one-parameter families (for certain combinations of $s$ and $m$) of profiles for~\eqref{eq:fpme1} or~\eqref{eq:fpme3}.

The special types of Barenblatt profiles sought here are proportional to $(R^2+|y|^2)^{-q}$ or $(R^2-|y|^2)_+^{q}$, for some $R>0$ and $q>0$. This is motivated from the Barenblatt profiles of the classic PME $u_t=
\Delta u^m$, which take the form of $(R^2+|y|^2)^{-1/(1-m)}$ for $m \in (\frac{N-2}{N},1)$ or $(R^2-|y|^2)_+^{1/(m-1)}$ for $m >1$. The main result is summarized as follows. 

\begin{subequations}
For~\eqref{eq:fpme1}, three families of explicit self-similar solutions of the form $(R^2+|y|^2)^{-q}$ are found
for $s\in (0,1)$: 
\begin{enumerate}
 \item when $m=\frac{N+2-2s}{N+2s} >m_c:=\frac{N-2s}{N}$, 
\begin{equation}\label{eq:fpme1MC}
 u(x,t) = \lambda t^{-N\beta}\big(R^2+|xt^{-\beta}|^2\big)^{-s-\frac{N}{2}},
\qquad \beta = \frac{1}{N(m-1)+2s};
\end{equation}
\item when $m = \frac{N-2s}{N+2s}<m_c$, 
\begin{equation}\label{eq:fpme1FTE}
 u(x,t) = \lambda (T-t)^{\frac{N+2s}{4s}}\big(R^2+|x|^2\big)^{-\frac{N}{2}-s};
\end{equation}
\item when $m = \frac{N-2s}{N+2s-2}$,
\begin{equation}\label{eq:fpme1IM}
 u(x,t) = \lambda t^{-\frac{N+2s-2}{2(1-s)}}
\big(R^2+|xt^{-\frac{1}{2(1-s)}}|^2\big)^{-\frac{N}{2}-s+1}.
\end{equation}
\end{enumerate}
\end{subequations}

For~\eqref{eq:fpme3}, only one family of self-similar of explicit self-similar solutions of the form $(R^2+|y|^2)^{-q}$ is found 
for $s\in (0,1)$:  when $m = \frac{N+6s-2}{N+2s}$, 
\begin{equation}\label{eq:fpme3MC}
 u(x,t) = \lambda t^{-N\beta} \big(R^2+|xt^{-\beta}|^2\big)^{-\frac{N}{2}-s},\qquad 
\beta = \frac{1}{N(m-1)+2-2s}.
\end{equation}

To derive these Barenblatt profiles, we need some preliminary results related to hypergeometric functions and their fractional Laplacians, given in Section~\ref{sec:wsint}.
The mass conserving Barenblatt profiles~\eqref{eq:fpme1MC} for~\eqref{eq:fpme1}
 are constructed in Section~\ref{sec:fpme1}, followed by 
mass conserving Barenblatt profiles~\eqref{eq:fpme3MC} for~\eqref{eq:fpme3}
in Section~\ref{sec:fpme3}.
The more complicated Barenblatt profiles~\eqref{eq:fpme1FTE} and ~\eqref{eq:fpme1IM} for~\eqref{eq:fpme1} with second-kind self-similarity are derived in Section~\ref{sec:Afpme1}.

%%%%%%%%%%%%%%%%%%%%%%%%%%%%%%%%%%%%%%%%%%%%%%%%%%%%%%%%%%%%%%%%%%%%%%%%%%%%%%%%%
%%%%%%%%%%%%%%%%%%%%%%%%%%%%%%%%%%%%%%%%%%%%%%%%%%%%%%%%%%%%%%%%%%%%%%%%%%%%%%%%%
\section{Fractional Laplacians of the Barenblatt profiles and other identities}
\label{sec:wsint}

In the search of Barenblatt profiles of the form $\Phi(y)=(R^2-|y|^2)_+^q$ or $\Phi(y)=(R^2+|y|^2)^{-q}$, the explicit expressions for the fractional Laplacian of $\Phi(y)$ are derived using Fourier transform. Certain special functions enter during various stages of the derivation, and therefore their definitions with related properties are introduced here. Most of the properties used here can be consulted from standard textbooks on special functions~\cite{MR1688958}.

Bessel-type special functions appear in the Fourier transform of $\Phi(y)$.
The \emph{Bessel functions of the first kind} $J_\nu(x)$ is the solution of the Bessel differential equation
\begin{equation*}\label{eq:besselde}
    x^2\frac{d^2z}{dx^2}+x\frac{dz}{dx}+(x^2-\nu^2)z=0,
\end{equation*}
that is finite at the origin for positive $\nu$.  The \emph{modified Bessel function of the second kind} $K_\nu(x)$ is the
exponentially decaying solution of the modified Bessel differential
equation
\begin{equation*}\label{eq:modifiedbesselde}
    x^2\frac{d^2z}{dx^2}+x\frac{dz}{dx}-(x^2+\nu^2)z=0.
\end{equation*}
In fact, besides the definitions, the only property we use below is $K_{-\nu}(x) = K_\nu(x)$.

The (Gauss) \emph{hypergeometric function} appears in the fractional 
Laplacian of $\Phi(y)$, which is a solution of Euler's
hypergeometric differential equation
\[
    x(1-x)\frac{d^2z}{dx^2}
    +\big[c-(a+b+1)x\big]\frac{dz}{dx}-abz=0,
\]
for any complex number $a,b$ and $c$. It is often represented more 
conveniently as a power series
\begin{equation}\label{eq:hg}
    { }_2F_1(a,b;c;x) = \sum_{n=0}^\infty
    \frac{(a)_n(b)_n}{(c)_n n!}x^n,\qquad |x|<1,
\end{equation}
where $(a)_n = \Gamma(a+n)\big/\Gamma(a)$ is the Pochhammer symbol and 
$\Gamma(x)$ is the Euler Gamma function.
If $c$ is a non-positive integer, ${ }_2F_1(a,b;c;x)$ becomes a
polynomial of degree $-c$ in $x$. From the series expansion~\eqref{eq:hg}, it is obvious that 
${ }_2F_1(a,b;c;x)={ }_2F_1(b,a;c;x)$ and 
\begin{equation}\label{eq:hgdev}
\frac{d}{dx} { }_2F_1(a,b;c;x) = \frac{ab}{c}\ {
}_2F_1(a+1,b+1;c+1;x).
\end{equation}
These two simple properties still hold on the complex plane, by 
analytical continuation.

The hypergeometric function is prevalent in mathematical physics
because it represents many other common yet important special functions and
it emerges also in many special integrals. 
In fact, the candidate Barenblatt profiles $(R^2-|y|^2)_+^q$ or $(R^2+|y|^2)^{-q}$ 
are also special hypergeometric functions, i.e.,
\begin{align}\label{eq:spehg}
    (R^2-|y|^2)_+^q &= R^{2q}{ }_2F_1(-q,c;c;|y|^2/R^2),\\
    (R^2+|y|^2)^{-q} &= R^{-2q}{ }_2F_1(q,c;c;-|y|^2/R^2),
\end{align}
for any complex number $c$. In this paper, we will always choose
$c$ to be $N/2$, half of the space dimension, to match the parameters in the Fourier transforms of $\Phi(y)$. 

For the explicit Barenblatt profiles of~\eqref{eq:fpme2} found in~\cite{biler2011barenblatt,biler2013nonlocal},
the key formula is  the \emph{Weber-Schafheitlin integral}~\cite[p. 401-403]{MR0010746}
\begin{multline*}
\int_0^\infty \eta^{-\rho}J_\mu(\eta a)J_{\nu}(\eta b)d\eta \cr =
\frac{b^{\nu}a^{\rho-\nu-1}\Gamma\left(\frac{\nu-\rho+\mu+1}{2}\right)}
{2^\rho \Gamma(\nu+1)\Gamma\left(\frac{1+\mu+\rho-\nu}{2}\right)} 
{}_2F_1\left(
\frac{\nu-\rho+\mu+1}{2},\frac{\nu-\rho-\mu+1}{2};\nu+1;\frac{b^2}{a^2}
\right),
\end{multline*}
with $\nu+\mu-\rho+1>0$, $\rho >-1$ and $0<b\leq a$. It
enables the authors to derive explicitly the (inverse) fractional Laplacian of $(R^2-|y|^2)_+^{q}$ for any $q>0$, $s \in (0,1)$, i.e.,
\begin{align}\label{eq:flmp}
&\quad (-\Delta)^{-s} \big( (R^2-|y|^2)_+^{q}\big)\cr
   &= \begin{cases}       
    C_{q,s,N} R^{2q+2s}\ {
    }_2F_1\left(\frac{N}{2}-s,-q-s;\frac{N}{2};|y|^2/R^2\right), \qquad 
    &|y|\leq R,\cr
    \tilde{C}_{q,s,N} R^{N+2q}|y|^{2s-N}{ }_2F_1\left(
    \frac{N}{2}-s,1-s;\frac{N}{2}+q+1;R^2/|y|^2\right),
    &|y|\geq R,
\end{cases}
\end{align}
with
\[
    C_{q,s,N} = \frac{2^{-2s}\Gamma(q+1)\Gamma(N/2-s)}{
    \Gamma(N/2)\Gamma(q+s+1)},\quad
    \tilde{C}_{q,s,N} = 
    \frac{2^{-2s}\Gamma(q+1)\Gamma(N/2-s)}{\Gamma(s)
    \Gamma(N/2+q+1)}.
\]

In this paper, we obtain explicit expressions for the fractional Laplacians  of $(R^2+|y|^2)^{-q}$, using the closely related
 \emph{modified Weber-Schafheitlin integral}~\cite[p. 410]{MR0010746}, 
which reads
\begin{multline}\label{eq:mdws}
\int_0^\infty \eta^{-\rho} K_\mu(\eta a) J_\nu(\eta b) d\eta \cr
= \frac{b^{\nu}a^{\rho-\nu-1}\Gamma\left(\frac{\nu-\rho+\mu+1}{2}\right)
\Gamma\left(\frac{\nu-\rho-\mu+1}{2}\right)}{2^{\rho+1}\Gamma(\nu+1)} 
{}_2F_1\left(
\frac{\nu-\rho+\mu+1}{2},\frac{\nu-\rho-\mu+1}{2};\nu+1;-\frac{b^2}{a^2}
\right),
\end{multline}
with $|\mu|< \nu-\rho+1$ and $a>0$. 

Already observed in~\cite{biler2011barenblatt,biler2013nonlocal}, these special Weber-Schafheitlin integrals are connected to the fractional Laplacians of $(R^2-|y|^2)_+^q$ or $(R^2+|y|^2)^{-q}$ by the fact that the Fourier transform of ${}_2F_1\big(a,b;\frac{N}{2};\pm |y|^2\big)$ are  $J_\nu(|\xi|)$ or $K_\nu(|\xi|)$, multiplied with a power of
$|\xi|$. The Barenblatt profiles $\Phi(y)=(R^2+|y|^2)^{-q}$ we are interested in this paper can be written as $R^{-2q}{}_2F_1\big(q,\frac{N}{2};\frac{N}{2};-\frac{|y|^2}{R^2}\big)$,
suggesting the choices of parameters $\nu = \frac{N}{2}-1$,
$a=R$ and $b = |y|$ in~\eqref{eq:mdws} while the rest two parameters $\mu$
and $\rho$ are chosen according to other parameters like $q$ and $N$. Comparing the expressions of the inverse Fourier transform of
radial functions given by~\eqref{eq:radialift} in Appendix~\ref{app:ft}, we conclude the following 
Fourier transform pair
\begin{equation*}
{ }_2F_1\left(\frac{N}{4}+\frac{\mu-\rho}{2},\frac{N}{4}
-\frac{\mu+\rho}{2};\frac{N}{2};-\frac{|y|^2}{R^2}\right), \quad 
\frac{2^{\rho+1}(2\pi)^{\frac{N}{2}}\Gamma(\frac{N}{2}) R^{\frac{N}{2}-\rho}}
{\Gamma(\frac{N}{4}+\frac{\mu-\rho}{2})
\Gamma(\frac{N}{4}-\frac{\mu+\rho}{2})}|\xi|^{-\rho-\frac{N}{2}}
K_\mu(|\xi| R).
\end{equation*}
This Fourier pair once again implies the following relation (with some restrictions on the parameters $\rho$, $\mu$ and $s$) for the 
fractional Laplacian of general hypergeometric functions
\begin{multline*}\label{eq:flhypergeom}
(-\Delta)^{s}\left[ 
{ }_2F_1\left(\frac{N}{4}+\frac{\mu-\rho}{2},\frac{N}{4}
-\frac{\mu+\rho}{2};\frac{N}{2};-\frac{|y|^2}{R^2}\right)
\right] \cr = 
2^{2s}R^{-2s} \frac{\Gamma(\frac{N}{4}+\frac{\mu-\rho}{2}+s)\Gamma(\frac{N}{4}
-\frac{\mu+\rho}{2}+s)}{\Gamma(\frac{N}{4}
+\frac{\mu-\rho}{2})\Gamma(\frac{N}{4}
-\frac{\mu+\rho}{2})} 
{ }_2F_1\left(\frac{N}{4}+\frac{\mu-\rho}{2}+s,\frac{N}{4}
-\frac{\mu+\rho}{2}+s;\frac{N}{2};-\frac{|y|^2}{R^2}\right).
\end{multline*}
In particular, when $\rho=-q$ and $\mu=\frac{N}{2}-q$, we get 
\begin{align}
 (-\Delta)^s (R^2+|y|^2)^{-q} &= 
R^{-2q} (-\Delta)^s \left[{}_2F_1\Big(q,\frac{N}{2};\frac{N}{2};-\frac{|y|^2}{R^2}\Big)\right] \cr
&=2^{2s}R^{-2s-2q}
\frac{\Gamma(q+s)\Gamma(\frac{N}{2}+s)}
{\Gamma(q)\Gamma(\frac{N}{2})}
{}_2F_1\Big(q+s,\frac{N}{2}+s;\frac{N}{2};-\frac{|y|^2}{R^2}\Big). 
\end{align}
For the explicit Barenblatt profiles we find for~\eqref{eq:fpme1}
and~\eqref{eq:fpme3} below, only two simple cases of~\eqref{eq:flhypergeom} are needed, which are collected here:
\begin{enumerate}[(i)]
\item when $q=\frac{N}{2}+1-s$, 
\begin{multline}\label{eq:expflfpme1}
    (-\Delta)^{s} (R^2+|y|^2)^{-\frac{N}{2}-1+s} \\
= 2^{2s-1}NR^{-N-2}\frac{\Gamma(\frac{N}{2}+s)}
        {\Gamma(\frac{N}{2}+1-s)}\ 
{}_2F_1\left(\frac{N}{2}+1,\frac{N}{2}+s;\frac{N}{2};-\frac{|y|^2}{R^2}\right); \qquad
\end{multline}
\item when $q=\frac{N}{2}-s$, 
\begin{equation}\label{eq:flpower}
    (-\Delta)^{s}  (R^2+|y|^2)^{-\frac{N}{2}+s}
=2^{2s}R^{2s}\frac{\Gamma(\frac{N}{2}+s)}
{\Gamma(\frac{N}{2}-s)}
(R^2+|y|^2)^{-\frac{N}{2}-s}.
\end{equation}
\end{enumerate}

\begin{remark}
Since there is no restriction on the sign of $s$, the inverse fractional 
Laplacian $(-\Delta)^{-s}$ of above functions  can be 
obtained by changing $s$ to $-s$. 
\end{remark}

In the next three sections, we search for Barenblatt profiles $\Phi(y)$ of the form
$\lambda(R^2+|y|^2)^{-q}$ or $\lambda(R^2-|y|^2)_+^{q}$, by looking at the local power series
expansion of the governing equation for $\Phi(y)$ at the origin. Moreover, for mass conserving 
self-similar solutions in the next two sections, the governing equation 
can be simplified to an identity involving two Gauss hypergeometric functions. 
The corresponding profiles are obtained using the following lemma, which is proved easily also using a power series expansion at the origin.
\begin{lemma}\label{lem:hypereq}
If the non-constant hypergeometric functions ${}_2F_1(a_1,b_1;c;x)$ and ${}_2F_1(a_2,b_2;c;x)$ 
are identical for $|x|<1$, then either $a_1=a_2,b_1=b_2$ or $a_1=b_2,b_1=a_2$.
\end{lemma}

%%%%%%%%%%%%%%%%%%%%%%%%%%%%%%%%%%%%%%%%%%%%%%%%%%%%%%%%%%%%%%%%%%%%%%%%%%%%%%%%%

%%%%%%%%%%%%%%%%%%%%%%%%%%%%%%%%%%%%%%%%%%%%%%%%%%%%%%%%%%%%%%%%%%%%%%%%%%%%%%%%%
%%%%%%%%%%%%%%%%%%%%%%%%%%%%%%%%%%%%%%%%%%%%%%%%%%%%%%%%%%%%%%%%%%%%%%%%%%%%%%%%%
\section{Mass conserving Barenblatt profiles for \texorpdfstring{$u_t +(-\Delta)^{s}u^m=0$}{mc}} \label{sec:fpme1}

If $u(x,t)$ is a solution of~\eqref{eq:fpme1}, so is 
$T_\lambda u(x,t) = \lambda^{N\beta}u(\lambda^\beta x, \lambda t)$
with 
\begin{equation}
\beta = \frac{1}{N(m-1)+2s}.
\end{equation}
This implies self-similar solutions of the form $ u(x,t) = t^{-N\beta}\Phi(y)$ with
$y = xt^{-\beta}$, where the Barenblatt profile $\Phi$ satisfies the equation
\begin{equation}\label{eq:Bfpme1}
    (-\Delta)^s \Phi^m = \beta \nabla\cdot(y \Phi).
\end{equation}
The basic existence, uniqueness and many properties of $\Phi(y)$ are already 
established by V\'{a}zquez~\cite{vazquez2012barenblatt}, without any explicit 
expressions of $\Phi(y)$ (except the linear case $m=1$ and $s=1/2$). Since the solutions~\eqref{eq:fpme1} become positive instantaneously~\cite{de2012general}, we do not expect Barenblatt profiles of the form $\Phi(y)=\lambda(R^2-|y|^2)_+^q$ and hence concentrate on $\Phi(y)=\lambda(R^2+|y|^2)^{-q}$ only. In fact, we have the following theorem.
\begin{theorem}
 For every $s\in (0,1)$, equation~\eqref{eq:fpme1} admits a self-similar solution 
$u(x,t)=t^{-N\beta}\Phi(xt^{-\beta})$ with the special profile $\Phi(y)=\lambda (R^2+|y|^2)^{-q}$
($q>0$) and $\beta = \frac{1}{N(m-1)+2s}$ only when $m=\frac{N+2-2s}{N+2s}$. The corresponding self-similar solution
\[
  u(x,t) = \lambda t^{-N\beta}\big(R^2+|xt^{-\beta}|^2\big)^{-s-\frac{N}{2}}
\]
is a classical solution on $(0,\infty)\times \mathbb{R}^N$ with $u(x,t)\to M\delta(x)$ as $t\to 0$ for some $M>0$.
\end{theorem}

To derive the Barenblatt profile,  replacing $q$
with $mq$ in~\eqref{eq:flhypergeom}, 
\[
 (-\Delta)^s\Phi(y)^m = 
\lambda^m 2^{2s}R^{-2s-2mq}\frac{\Gamma(mq+s)\Gamma(\frac{N}{2}+s)}{\Gamma(mq)\Gamma(\frac{N}{2})}
{}_2F_1\big(mq+s,\frac{N}{2}+s;\frac{N}{2};-\frac{|y|^2}{R^2}\big).
\]
On the other hand, a simple calculation yields
\[
 \nabla\cdot \big(y\Phi(y)\big)=\lambda N R^{-2q}{}_2F_1\big(q,\frac{N}{2}+1;\frac{N}{2};-\frac{|y|^2}{R^2}\big).
\]
As a result, the governing equation~\eqref{eq:Bfpme1} reduces to the identity
\begin{equation}\label{eq:fpme1hypereq}
 {}_2F_1\big(mq+s,\frac{N}{2}+s;\frac{N}{2};-\frac{|y|^2}{R^2}\big)={}_2F_1\big(q,\frac{N}{2}+1;\frac{N}{2};-\frac{|y|^2}{R^2}\big)
\end{equation}
and the algebraic equation
\begin{equation}\label{eq:fpme1algeq}
 \lambda^m 2^{2s}R^{-2s-2mq}\frac{\Gamma(mq+s)\Gamma(\frac{N}{2}+s)}{\Gamma(mq)\Gamma(\frac{N}{2})} = \beta \lambda N R^{-2q}.
\end{equation}
Since $\frac{N}{2}+s\neq \frac{N}{2}+1$ in~\eqref{eq:fpme1hypereq}, Lemma~\ref{lem:hypereq} implies that
\[
 mq+s=\frac{N}{2}+1,\qquad \frac{N}{2}+s=q,
\]
or 
\begin{equation}\label{eq:mqfpme1}
 m = \frac{N+2-2s}{N+2s},\qquad q=\frac{N}{2}+s.
\end{equation}
Consequently, the algebraic identity~\eqref{eq:fpme1algeq} can be simplified as
\begin{equation}\label{eq:algfpme1}
\lambda^{1-m}R^{2-2s}\beta = 2^{2s-1}\frac{\Gamma(\frac{N}{2}+s)}
{\Gamma(\frac{N}{2}+1-s)}.
\end{equation}
Together with the total mass condition
\begin{equation}\label{eq:masscond}
M = \int_{\mathbb{R}^N} \Phi(y)dy = 
\lambda \pi^{\frac{N}{2}}R^{-2s}\frac{\Gamma(s)}{\Gamma(\frac{N}{2}+s)},
\end{equation}
the two free parameters $\lambda$ and $R$ are determined uniquely.

\begin{remark}
The special case $m=1$ and $s=1/2$ is well-known,  and the corresponding Barenblatt profile is the  Poisson kernel. The new solutions above can be viewed as a continuous branch from the point
$s=1/2$ to the whole interval $s \in (0,1)$.
\end{remark}

\begin{remark}
These Barenblatt profiles are obtained for $m=\frac{N+2-2s}{N+2s}>m_c := \frac{(N-2s)_+}{N}$, and 
have the solutions $u(\cdot,t)\in L^1(\mathbb{R}^N)$ for any $t>0$. 
The general functional framework of existence and uniqueness developed in~\cite{de2012general} applies here. Moreover, the optimal decay rate $O(|y|^{-N-2s})$ of general Barenblatt profiles governed by~\eqref{eq:Bfpme1} 
is proved in~\cite{bonforte2012quantitative,vazquez2012barenblatt} for $m>m_1:=\frac{N}{N+2s}$,
which is also verified in above special cases since $m=\frac{N+2-2s}{N+2s}>m_1$.
\end{remark}

%%%%%%%%%%%%%%%%%%%%%%%%%%%%%%%%%%%%%%%%%%%%%%%%%%%%%%%%%%%%%%%%%%%%%%%%%%%%%%%%%

%%%%%%%%%%%%%%%%%%%%%%%%%%%%%%%%%%%%%%%%%%%%%%%%%%%%%%%%%%%%%%%%%%%%%%%%%%%%%%%%%
%%%%%%%%%%%%%%%%%%%%%%%%%%%%%%%%%%%%%%%%%%%%%%%%%%%%%%%%%%%%%%%%%%%%%%%%%%%%%%%%%
\section{Mass conserving Barenblatt profiles \texorpdfstring{$u_t=\nabla\cdot(u^{m-1}\nabla (-\Delta)^{-s}u)$}{nlpre}}\label{sec:fpme3}

Since solutions of~\eqref{eq:fpme3} could have either finite (for $m\geq 2$) or infinite speed of propagation (for $1<m<2$) as shown in~\cite{prop}, 
Barenblatt profiles of both forms $\lambda(R^2+|y|^2)^{-q}$ and $\lambda(R^2-|y|^2)_+^q$ are sought in this section.

If $u(x,t)$ is a solution of~\eqref{eq:fpme3}, so is $T_\lambda u(x,t) = \lambda^{N\beta}u(\lambda^\beta x, \lambda t)$
with 
\begin{equation}
\beta = \frac{1}{N(m-1)+2-2s}.
\end{equation}
This implies self-similar solutions of the form $ u(x,t) = t^{-N\beta}\Phi(y)$ with
$y = xt^{-\beta}$, where the Barenblatt profile $\Phi$ satisfies 
\begin{equation}\label{eq:Bfpme3F}
  \nabla\cdot\big(\Phi^{m-1}\nabla(-\Delta)^{-s}\Phi\big)+ \beta \nabla\cdot\big(y \Phi\big) = 0.
\end{equation}
Since the special case $m=2$ is already covered in~\cite{biler2011barenblatt,biler2013nonlocal}, 
we only consider the case $m\neq 2$ below.
\begin{theorem}
 If $m\neq 2$, for every $s\in (0,1)$, equation~\eqref{eq:fpme3} admits a self-similar solution $u(x,t)=
t^{-N\beta}\Phi(xt^{-\beta})$ with the special profile $\Phi(y)=\lambda (R^2+|y|^2)^{-q}$ ($q>0$)
only when $m = \frac{N+6s-2}{N+2s}$. The corresponding self-similar solution
\[
 u(x,t) = \lambda t^{-N\beta}\big(R^2+|xt^{-\beta}|^2\big)^{-s-\frac{N}{2}}
\]
is a classical solution on $(0,\infty)\times \mathbb{R}^N$ with $u(x,t) \to M\delta(x)$
as $t\to 0$ for some $M>0$. Furthermore, equation~\eqref{eq:fpme3} does not admit
any self-similar solution $u(x,t)=t^{-N\beta}\Phi(xt^{-\beta})$ with the special profile $\Phi(y)=\lambda (R^2-|y|^2)_+^{q}$.
\end{theorem}

To facilitate the calculation,  the governing equation~\eqref{eq:Bfpme3F} can be integrated once and then simplified as
\begin{equation}\label{eq:Bfpme3}
 \nabla(-\Delta)^{-s}\Phi+\beta y\Phi^{2-m}=0,
\end{equation}
whenever $\Phi\neq 0$.

\subsection{Barenblatt profiles of the form $\Phi(y)=\lambda (R^2+|y|^2)^{-q}$}

In this case we get by~\eqref{eq:flhypergeom}
\[
(-\Delta)^{-s}\Phi(y)=\lambda 2^{-2s}R^{2s-2q}\frac{\Gamma(q-s)\Gamma(\frac{N}{2}-s)}{\Gamma(q)\Gamma(\frac{N}{2})}
{}_2F_1\big(q-s,\frac{N}{2}-s;\frac{N}{2};-\frac{|y|^2}{R^2}) 
\]
and consequently $\nabla (-\Delta)^{-s}\Phi(y)$ becomes
\[
-\lambda 2^{1-2s}R^{2s-2q-2}\frac{\Gamma(q-s+1)\Gamma(\frac{N}{2}-s+1)}{\Gamma(q)\Gamma(\frac{N}{2}+1)}
y{}_2F_1\big(q-s+1,\frac{N}{2}-s+1;\frac{N}{2}+1;-\frac{|y|^2}{R^2}).
\]
On the other hand, $\Phi^{2-m}$ can be written as
\[
 \lambda^{2-m}(R^2+|y|^2)^{-q(2-m)}
=\lambda^{2-m}R^{-2q(2-m)}{}_2F_1\big(q(2-m),\frac{N}{2}+1;\frac{N}{2}+1;-\frac{|y|^2}{R^2}\big).
\]
Therefore, the simplified governing equation~\eqref{eq:Bfpme3} reduces to the identity
\begin{equation}\label{eq:fpme3IF}
 {}_2F_1\big(q-s+1,\frac{N}{2}-s+1;\frac{N}{2}+1;-\frac{|y|^2}{R^2})=
{}_2F_1\big(q(2-m),\frac{N}{2}+1;\frac{N}{2}+1;-\frac{|y|^2}{R^2}\big)
\end{equation}
and the algebraic equation
\begin{equation}\label{eq:fpme3AE}
 -\lambda 2^{1-2s}R^{2s-2q-2}\frac{\Gamma(q-s+1)\Gamma(\frac{N}{2}-s+1)}{\Gamma(q)\Gamma(\frac{N}{2}+1)}
+\beta\lambda^{2-m}R^{-2q(2-m)}=0.
\end{equation}
Since $\frac{N}{2}-s+1\neq \frac{N}{2}+1$,~\eqref{eq:fpme3IF} holds if and only if 
\[
q-s+1=\frac{N}{2}+1,\qquad 
\frac{N}{2}-s+1=q(2-m) 
\]
or
\begin{equation}
q=\frac{N}{2}+s,\qquad 
m = \frac{N+6s-2}{N+2s}.
\end{equation}
As a result,~\eqref{eq:fpme3AE} can be simplified as 
\[
 \lambda^{1-m}R^{2s}\beta = 2^{1-2s}\frac{\Gamma(\frac{N}{2}-s+1)}{\Gamma(\frac{N}{2}+s)},
\]
which determines $\lambda$ and $R$ uniquely, together with~\eqref{eq:masscond} for the total mass.

\subsection{Barenblatt profiles of the form $\Phi(y)=\lambda (R^2-|y|^2)_+^{q}$}
In this case, using~\eqref{eq:flmp},  $\nabla(-\Delta)^{-s}\Phi(y)$ for $|y|<R$ can be written as
\[
-2\lambda \frac{C_{q,s,N}(N-2s)(q+s)}{N}R^{2q+2s-2}y{}_2F_1\left(
\frac{N}{2}-s+1,-q-s+1;\frac{N}{2}+1;\frac{|y|^2}{R^2}
\right).
\]
On the other hand,
\[
 \beta y\Phi(y)^{2-m}=
\beta y (R^2-|y|^2)_+^{(2-m)q}
=\beta yR^{2(2-m)q}{}_2F_1\left(-(2-m)q,\frac{N}{2}+1;\frac{N}{2}+1;\frac{|y|^2}{R^2}\right).
\]
Therefore, the simplified governing equation~\eqref{eq:Bfpme3}
is satisfied only if
\[
 {}_2F_1\left(
\frac{N}{2}-s+1,-q-s+1;\frac{N}{2}+1;\frac{|y|^2}{R^2}
\right)
={}_2F_1\left(-(2-m)q,\frac{N}{2}+1;\frac{N}{2}+1;\frac{|y|^2}{R^2}\right).
\]
Since $m\neq 2$, the hypergeometric function on the right hand side is non-constant. By Lemma~\ref{lem:hypereq}, we must have 
\[
 -q-s+1=\frac{N}{2}+1,\quad 
\frac{N}{2}-s+1=-(2-m)q.
\]
Since both $q$ and $s$ are positive, the first equation is invalid and there is no Barenblatt profiles of these equations. Therefore, there is no non-trivial Barenblatt profiles of the type
$\lambda(R^2-|x|^2)_+^q$ when $m>2$, despite the existence of solutions
propagating with finite speed in one dimension~\cite{prop}.

%%%%%%%%%%%%%%%%%%%%%%%%%%%%%%%%%%%%%%%%%%%%%%%%%%%%%%%%%%%%%%%%%%%%%%%%%%%%%%%%%
%%%%%%%%%%%%%%%%%%%%%%%%%%%%%%%%%%%%%%%%%%%%%%%%%%%%%%%%%%%%%%%%%%%%%%%%%%%%%%%%%
\section{Second-kind Barenblatt profiles  for \texorpdfstring{$u_t + (-\Delta)^s
u^m=0$}{2nd}}    \label{sec:Afpme1}

In the previous two sections, explicit self-similar solutions $u(x,t) = t^{-\alpha}\Phi(xt^{-\beta})$ are sought with
the \emph{a priori} condition $\alpha = N\beta$,
reflecting the mass conservation of these special solutions. However, this condition may break down, leading to the 
concept of self-similar solutions of the \emph{second kind}~\cite{MR1426127}. 
For~\eqref{eq:fpme1}, these anomalous self-similar solutions could appear in two situations.
In the fast diffusion regime $m<(N-2s)_+/N$, it is known that the mass escapes to infinity and the solution becomes identically zero at some finite time $T$~\cite{de2012general}. Here the self-similar solution, if it exists, takes the form
\begin{equation}\label{eq:fpme1F}
u(x,t) = (T-t)^{\alpha}\Phi\big( x(T-t)^\beta\big),
\end{equation}
with the restriction $\alpha > N\beta$.
On the other hand, the solution may have infinite mass, and hence it does not make 
any sense to require the solution to "conserve" the total mass. Here the self-similar solution takes the form
\begin{equation}\label{eq:fpme1M}
u(x,t) = t^{-\alpha}\Phi\big( xt^{-\beta}\big),
\end{equation}
where $\Phi(y)$ decays slower than $|y|^{-N}$ as $|y|\to \infty$ and the relation between $\alpha$ and $\beta$
cannot be determined \emph{a priori}. Since the Barenblatt profiles $\Phi$ for both~\eqref{eq:fpme1F}
and~\eqref{eq:fpme1M} satisfy the same equation
\begin{equation}\label{eq:2ndBfpme1}
(-\Delta)^s \Phi^m - \alpha \Phi - \beta y\cdot \nabla \Phi=0,
\end{equation}
we treat them at the same time below. Notice that there is only one condition on the scaling exponents $\alpha$ and
$\beta$, i.e., $\alpha(m-1)+2s\beta=-1$ for~\eqref{eq:fpme1F} or $\alpha(m-1)+2s\beta=1$ for~\eqref{eq:fpme1M}, which
is not enough to determine $\alpha$ and $\beta$ explicitly as in the previous two sections.

\begin{theorem}\label{thm:2nd}
For every $s\in (0,1)$, equation~\eqref{eq:fpme1} admits two self-similar solutions of the 
\emph{second kind} with profile $\Phi(y)=\lambda (R^2+|y|^2)^{-q}$: 
\begin{enumerate}[(a)]
 \item when $m = \frac{N-2s}{N+2s}$, the self-similar solution
\begin{equation}\label{eq:2ndFTE}
 u(x,t) = \lambda (T-t)^{\frac{N+2s}{4s}}\big(R^2+|x|^2\big)^{-\frac{N}{2}-s}
\end{equation}
is a classical solution on $[0,T)\times \mathbb{R}^N$ and vanishes at finite time $T>0$.

\item when $m=\frac{N-2s}{N+2s-2}$, the self-similar solution
\begin{equation}\label{eq:2ndVM}
 u(x,t) = \lambda t^{-\frac{N+2s-2}{2(1-s)}}
\big(R^2+|xt^{-\frac{1}{2(1-s)}}|^2\big)^{-\frac{N}{2}-s+1}
\end{equation}
is a classical solution on $(0,\infty)\times \mathbb{R}^N$ and has infinite mass at any $t>0$.
\end{enumerate}
\end{theorem}

Because $\alpha$ is different from $N\beta$ for the second kind self-similarity, the three
terms in the governing equation~\eqref{eq:2ndBfpme1} can not be simplified as an equation 
with two hypergeometric functions as in the previous two sections. Instead, we proceed in two steps. In the first step, we focus on the relation between the parameter $m$ and the exponent $q$ in the  rescaled profile $\Phi(y)=(1+|y|^2)^{-q}$ by a local series expansion for $r=|y|$.  
Using these explicit values of $m$ and $q$, we get the condition on $\lambda$ and $R$
in the general profile $\Phi(y)=\lambda (R^2+|y|^2)^{-q}$. In fact, the same two steps 
can be applied in Section~\ref{sec:fpme1}, to find the relation~\eqref{eq:mqfpme1} from the identity~\eqref{eq:fpme1hypereq} by power series expansions and the condition~\eqref{eq:algfpme1}
between $\lambda$ and $R$. 

When the simple, rescaled profile $\Phi(y)=(1+|y|^2)^{-q}$ is used, the governing equation~\eqref{eq:2ndBfpme1} should be rescaled too. The key observation is that, because the last two terms $\alpha\Phi$ and $\beta y\cdot \nabla\Phi$ have the same scaling factor, the relation
between $m$ and $q$ can be computed from $g(r)=0$, where $g(r)$ is a rescaled  version of~\eqref{eq:algfpme1} in the radial variable $r=|y|$, i.e., 
\begin{align}\label{eq:geq}
    g(r) = {}_2F_1\left(mq+s,\frac{N}{2}+s;\frac{N}{2};-r^2\right)
-(1+r^2)^{-q} - \tilde{\beta} r\frac{d}{dr}(1+r^2)^{-q},
\end{align}
for some $\tilde{\beta}$. Here the coefficient of $(-\Delta)^{s}\Phi$ or
$\alpha \Phi$ is scaled to unit, such that $g(0)=0$. The scale invariance of 
$y\cdot\nabla \Phi(y)/\Phi(y)$ implies that $\tilde{\beta}=\beta/\alpha$, which 
should be different from $1/N$ for the second kind self-similar solutions we are looking for here.
This scaling technique enables us to get the relation between $m$ and $q$, without worrying too much about the complicated constants or prefactors, while the remaining parameters in the Barenblatt profiles are then determined, using only the relatively simple identities~\eqref{eq:expflfpme1} or~\eqref{eq:flpower}. 

Finally, we can find the conditions that $g(r)$ vanishes identically from a 
power series expansion around the origin\footnote{A computer algebra system like MAPLE or MATHEMATICA is recommended to perform these symbolic calculations.}. that is 
$g(r) = g_0 + g_2r^2+g_4r^4+\cdots$.
Obviously $g_0$ vanishes. From 
\[
g_2 = {\frac {2\,\tilde{\beta}\,qN-Nmq-2\,mqs+Nq-Ns-2\,{s}^{2}}{N}}=0,
\]
we get
\[
m = {\frac {2\,\tilde{\beta}\,qN+Nq-Ns-2\,{s}^{2}}{q \left( N+2\,s \right) }}.
\]
Solving $q$ from 
{\small
\[
g_4={\frac { \left( 2\,{N}^{2}{\tilde{\beta}}^{2}q+4\,N{\tilde{\beta}}^{2}qs+4\,N{\tilde{\beta}}^{
2}q-{N}^{2}\tilde{\beta}-8\,\tilde{\beta}\,qs+4\,\tilde{\beta}\,{s}^{2}-2\,N\tilde{\beta}+Ns-4\,\tilde{\beta}
\,s-2\,qs+2\,{s}^{2} \right) q}{ \left( N+2 \right)  \left( N+2\,s
 \right) }}=0,
\]
}
to obtain (the other solution $q=0$ is irrelevant)
\[
q=\frac{1}{2}\,{\frac { \left( N+2\,s \right)  \left( N\tilde{\beta}-2\,\tilde{\beta}\,s+2\,
\tilde{\beta}-s \right) }{{N}^{2}{\tilde{\beta}}^{2}+2\,N{\tilde{\beta}}^{2}s+2\,N{\tilde{\beta}}^{2}-
4\,\tilde{\beta}\,s-s}}.
\]
Using the explicit expressions of $m$ and $q$, the coefficient $g_6$ can be simplified as
\begin{align*}
\frac{1}{3}
\frac{s(2\tilde{\beta}+1)(N+2s+2)(N+2s)(N\tilde{\beta}-2\tilde{\beta} s+2\tilde{\beta}-s)}{ \left( N+4 \right)  \left( {N}^{2}{\tilde{\beta}}^{2}+2\,N{\tilde{\beta}}^{2}s+2\,N{\tilde{\beta}}^{2}-4\,\tilde{\beta}\,s-s
 \right) ^{3}}.
\end{align*}
Here all the non-zero
factors in $g_6$ are isolated in the fractions, especially $N\tilde{\beta}-2\tilde{\beta} s+2\tilde{\beta}-s$
(otherwise $q=0$). We discuss the different cases for $g_6=0$ below, or all $\tilde{\beta}$
such that
\[
 \tilde{\beta} (N\tilde{\beta}-1)(N\tilde{\beta} -s)(N\tilde{\beta}+2\tilde{\beta} s-2\tilde{\beta}-1)=0.
\]

\begin{description}
\item[Case $\tilde{\beta} = 0$.] Then $m = \frac{N-2s}{N+2s}$, $q=\frac{N}{2}+s$ and $\alpha=\pm\frac{1}{1-m} =
    \pm\frac{N+2s}{4s}$. We have to choose $\alpha=\frac{N+2s}{4s}>0$, otherwise the corresponding self-similar
    solutions are growing in time, leading to the self-similar solution $u(x,t)=(T-t)^{\alpha}\Phi(x)$.
 The two constants $\lambda$ and $R$ in the Barenblatt profile $\Phi(y) = \lambda (R^2+|y|^2)^{-\frac{N}{2}-s}$ 
 are related by only one equation, the matching condition of coefficients from 
the identity~\eqref{eq:flpower}, i.e.,
\[
\lambda^{1-m}R^{2s} = \frac{\alpha}{2^{2s}}\frac{\Gamma(\frac{N}{2}-s)}{\Gamma(\frac{N}{2}+s)}.
\]
This gives the self-similar solution~\eqref{eq:2ndFTE} in Theorem~\ref{thm:2nd}, where 
$\lambda$ and $R$ can be determined uniquely by the initial mass
\[
M_0 = \int_{\mathbb{R}^N} u(x,0)dx= \lambda \pi^{\frac{N}{2}}T^{\frac{N+2s}{4s}}
\frac{\Gamma(s)}{\Gamma(\frac{N}{2}+s)}.                                                                                                                                            \]

\item[Case $N\tilde{\beta} -1=0$.] This implies that $\tilde{\beta} = 1/N=\beta/\alpha$ and 
it reduces the Barenblatt profiles considered in Section~\ref{sec:fpme1}.
\item[Case $N\tilde{\beta} -s =0$.] Then $q = -1<0$, leading to unacceptable solutions
growing at infinity.
\item[Case $N\tilde{\beta}+2\tilde{\beta} s-2\tilde{\beta}-1=0$ or $\tilde{\beta} = \frac{1}{N+2s-2}$.] The exponents $m$ and $q$ are simplified as
\[
m = \frac{N-2s}{N+2s-2},\quad 
q = \frac{N}{2}+s-1.
\]
The corresponding Barenblatt profiles $\Phi(y)=\lambda(R^2+|y|^2)^{-\frac{N}{2}-s+1}$ have infinite mass, as 
$q=\frac{N}{2}+s-1<\frac{N}{2}$. Since $m$ is strictly larger than $m_c=(N-2s)_+/N$,
the solutions do not vanish in finite, and  we expect the self-similar solutions~\eqref{eq:fpme1M} instead of~\eqref{eq:fpme1F}. This implies
$(m-1)\alpha+2s\beta=1$. Together with $\alpha = \beta/\tilde{\beta} = (N+2s-2)\beta$, we obtain
\[
\alpha = \frac{N+2s-2}{2(1-s)},\qquad 
\beta = \frac{1}{2(1-s)}.
\]
Finally, we find the relation $\lambda$ and $R$ in the profile $\Phi(y)=\lambda(R^2+|y|^2)^{-\frac{N}{2}-s+1}$. 
Since 
\[
(-\Delta)^s \Phi(y)^m = \lambda^m (-\Delta)^s (R^2+|y|^2)^{-\frac{N}{2}-s} 
=\lambda^m2^{2s}R^{2s}\frac{\Gamma(\frac{N}{2}+s)}{\Gamma(\frac{N}{2}-s)}(R^2+|y|^2)^{-\frac{N}{2}-s},
\]
and
\[
\alpha \Phi(y) + \beta y \cdot\nabla \Phi(y) = \alpha \lambda R^2 (R^2+|y|^2)^{-\frac{N}{2}-s},
\]
the equation~\eqref{eq:2ndBfpme1} for the profile is satisfied if
\[
\lambda^{1-m}R^{2-2s}=\frac{2^{2s}}{\alpha}\frac{\Gamma(\frac{N}{2}+s)}{\Gamma(\frac{N}{2}-s)}.
\]
This gives the self-similar solution~\eqref{eq:2ndVM} in Theorem~\ref{thm:2nd}, and in 
general $\lambda$ and $R$ can not be determined uniquely. 
\end{description}

The second kind self-similar solutions already appear in the literature in various contexts.
The finite-time extinction of solution for $m<m_c:=\frac{N-2s}{N}$ is already considered 
in~\cite{de2012general}, with estimates on the extinction time using functional 
inequalities~\cite{bonforte2012quantitative} or comparison in Marcinkiewicz norm~\cite{vazquez2013optimal}. For $m=\frac{N-2s}{N+2s}$, the self-similar solutions
constructed above is believed to better characterize the fine details
right before the extinction and provides a more accurate estimate on the extinction time 
for certain initial data. 

Solutions of~\eqref{eq:fpme1} with infinite mass do not fit into the general theoretical framework
for $L^1$ initial data developed in~\cite{MR2737788,de2012general} and have to be treated
in weighted space~\cite{bonforte2012quantitative}.  Therefore, in contrast to those first kind self-similar solutions starting with a Dirac delta initial condition, self-similar solutions
with singular initial data like $u(x,0)=|x|^{-N/p}$ for $p>\max(1,N(1-m)/2s)$ are shown to be
second kind, with conserved $L^p$ norm instead of $L^1$ norm (the mass). The solution~\eqref{eq:2ndVM} obtained above provides another explicit example of anomalous 
scaling for large data.

Finally, it should be noted another \emph{anomalous} self-similar solutions, so-called Very Singular Solutions (VSS), also constructed from separation of variables. when $0<m<m_c$, the solution
\[
u(x,t)=C(T-t)^{1/(1-m)}|x|^{-2s/(1-m)} 
\]
is used to estimate the finite extinction time~\cite{vazquez2013optimal}; when $\frac{N-2s}{N}:=m_c
<m<N/(N+2)$, the solution 
\[ 
u(x,t)=Ct^{1/(1-m)}|x|^{-2s/(1-m)}
\]
arises in the limit when the total mass of the first-kind Barenblatt profiles goes to 
infinity~\cite{vazquez2012barenblatt}. 
In the limit $R\to 0$,~\eqref{eq:2ndFTE} reduces to the former in the case $m=\frac{N-2s}{N+2s}$.
However, ~\eqref{eq:2ndVM} does not reduce to the latter because the range $m=\frac{N-2s}{N+2s-2}$
is not inside the interval $\big(m_c,N/(N+2)\big)$ of existence in general.
%%%%%%%%%%%%%%%%%%%%%%%%%%%%%%%%%%%%%%%%%%%%%%%%%%%%%%%%%%%%%%%%%%%%%%%%%%%%%%%%%
%%%%%%%%%%%%%%%%%%%%%%%%%%%%%%%%%%%%%%%%%%%%%%%%%%%%%%%%%%%%%%%%%%%%%%%%%%%%%%%%%
\section{Conclusion}

In this paper, several one-parameter families of explicit self-similar solutions are obtained for fractional porous medium equations~\eqref{eq:fpme1} and~\eqref{eq:fpme3}. The special forms of the Barenblatt profiles are motivated from the classic PMEs, and are determined from the matching conditions of certain hypergeometric functions or local power series expansions. 
These special scale invariant solutions can complement the qualitative and quantitative studies of the underlying equations with explicit examples, and provide immense intuition for further investigation. In addition, these exact solutions can also be used to test the accuracy and 
efficiency of numerical methods for equations with fractional operators.

By our construction, the explicit Barenblatt profiles are exhausted in the forms $\lambda(R^2+|y|^2)^{-q}$ or $\lambda(R^2-|y|^2)_+^q$ for the cases we sought. In contrast
to those of~\eqref{eq:fpme2} obtained for all $m$ and $s$ in~\cite{biler2011barenblatt,biler2013nonlocal}, these
explicit profiles for~\eqref{eq:fpme1} or~\eqref{eq:fpme3} exist only for certain combinations of $m$ and $s$. The profiles for general $m$ and $s$, whose existence 
may be relatively easy to prove as in~\cite{vazquez2012barenblatt}, are expected to have much more
complicated expressions (if they exist). The complexity can be observed from the explicit Barenblatt profiles of the fractional heat equation $u_t + (-\Delta)^s u=0$ via
Fourier transform. Therefore, it is interesting to see whether there are any explicit candidate profiles for the more general cases.

%%%%%%%%%%%%%%%%%%%%%%%%%%%%%%%%%%%%%%%%%%%%%%%%%%%%%%%%%%%%%%%%%%%%%%%%%%%%%%%%%

%%%%%%%%%%%%%%%%%%%%%%%%%%%%%%%%%%%%%%%%%%%%%%%%%%%%%%%%%%%%%%%%%%%%%%%%%%%%%%%%%%%%%%
\section*{Acknowledgements}
This work is supported by Engineering and Physical Sciences Research Council grant number EP/K008404/1.
The author would like to thank the hospitality of Professor Juan Luis V{\'a}zquez and
Universidad Aut{\'o}noma de Madrid where this work was initiated. The author also appreciates the 
anonymous referees for comments and suggestions to improve the paper.

%%%%%%%%%%%%%%%%%%%%%%%%%%%%%%%%%%%%%%%%%%%%%%%%%%%%%%%%%%%%%%%%%%%%%%%%%%%%%%%%%
%%%%%%%%%%%%%%%%%%%%%%%%%%%%%%%%%%%%%%%%%%%%%%%%%%%%%%%%%%%%%%%%%%%%%%%%%%%%%%%%%
\appendix

\section{Fourier transform of radial functions}
\label{app:ft}
The fractional Laplacian of Barenblatt profiles in this paper
is evaluated by Fourier 
transform and inverse Fourier transform. These transforms are defined as
\[
 \hat{u}(\xi) = \mathcal{F}[u](\xi)=
\int_{\mathbb{R}^N} u(x)e^{-i\xi\cdot x}dx,\qquad 
u(x) = \mathcal{F}^{-1}[u](x)=(2\pi)^{-N}
\int_{\mathbb{R}^N} \hat{u}(\xi)e^{i\xi\cdot x}dx.
\]
In particular, we need a few facts about the transforms of radially 
symmetry functions~\cite{MR2449250}.

Using explicit expression for the integration of $e^{i\omega\cdot x}$ over the unit sphere 
$\mathbb{S}^{N-1}$ in $\mathbb{R}^N$, i.e.,
\[
    \int_{\mathbb{S}^{N-1}} e^{i\omega\cdot x} d\omega 
    = (2\pi)^{\frac{N}{2}}|x|^{1-\frac{N}{2}}J_{\frac{N}{2}-1}(|x|),
\]
the Fourier transform of a radial function $u(|x|)$ becomes
\begin{equation}\label{eq:radialft}
\mathcal{F}[u](\xi) 
=(2\pi)^{\frac{N}{2}}|\xi|^{1-\frac{N}{2}}\int_0^\infty
    r^{\frac{N}{2}}J_{\frac{N}{2}-1}(r|\xi|)u(r)dr.
\end{equation}
Similarly the inverse Fourier transform of a radial function
$\hat{u}(|\xi|)$ becomes
\begin{equation}\label{eq:radialift}
\mathcal{F}^{-1}[\hat{u}](x) 
=(2\pi)^{-\frac{N}{2}}|x|^{1-\frac{N}{2}}\int_0^\infty
\eta^{\frac{N}{2}}J_{\frac{N}{2}-1}(\eta|x|)\hat{u}(\eta)d\eta.
\end{equation}

%\end{acknowledgements}

%\bibliography{biopme}
%\bibliographystyle{plain}

%\end{document}

%\affiliationone{% in this example, two authors share an institution
%   Yanghong Huang\\
%   Department of Mathematics\\
%   Imperial College London\\
%   London SW7 2AZ, UK
%   \email{yanghong.huang@imperial.ac.uk}}
% Important: Do not put any empty line here.
%\affiliationtwo{% in this example, one author has two addresses}
%   T. Hird\\
%   Previous postal address where
%     the research was performed and\\
%   Country
%   \email{hird@university.ac.uk}}
% Important: Do not put any empty line here.
% Use \affiliationthree{} for any address positioned under \affiliationone
% Use \affiliationfour{}  for any address positioned under \affiliationtwo
%\affiliationthree{~} %inserts a space to make this field empty
%\affiliationfour{%
%   Current address:\\
%   Present long-term address\\
%   Country
%   \email{t.hird@institution.edu}}
%
\end{document}